\documentclass[11pt]{article}

\usepackage{float}
\usepackage{multirow}
\usepackage{amsmath,amsfonts,amssymb,graphicx,amsthm,color,natbib}
\usepackage{authblk}
\allowdisplaybreaks
\usepackage[]{geometry}
\usepackage[normalem]{ulem}
\usepackage{hyperref}
\usepackage{booktabs}

\newtheorem{prop}{Proposition}%
\newtheorem{cor}{Corollary}

\begin{document}

	 \title{Censoring heavy-tail count distributions for parameter estimation with an application to stable distributions}

\author[1,2]{{Antonio} {Di Noia}}
\author[3]{{Marzia} {Marcheselli}}
\author[3]{{Caterina} {Pisani}}
\author[4]{{Luca} {Pratelli}}
\affil[1]{Seminar for Statistics, Department of Mathematics, ETH Zurich}
\affil[2]{Faculty of Economics, Università della Svizzera italiana}
\affil[3]{Department of Economics and Statistics, University of Siena}
\affil[4]{Italian Naval Academy}
\date{}
\setcounter{Maxaffil}{0}
\renewcommand\Affilfont{\itshape\small}
\maketitle
\let\thefootnote\relax\footnotetext{\emph{Email addresses:} antonio.dinoia(\sout{at})stat.math.ethz.ch (Antonio Di Noia), marzia.marcheselli(\sout{at})unisi.it (Marzia Marcheselli), caterina.pisani(\sout{at})unisi.it (Caterina Pisani), luca\_pratelli(\sout{at})marina.difesa.it (Luca Pratelli).}

\begin{abstract}
A new approach based on censoring and moment criterion is introduced for parameter estimation of count distributions when the probability generating function is available even though a closed form of the probability mass function and/or finite moments do not exist.
\end{abstract}


\section{Introduction}
\label{sec:intro}
Heavy-tailed count data naturally arise in many applied disciplines 
(see e.g., \citealp{ El-Shaarawi}, \citealp{Edwards}, \citealp{Sun}). A plethora of family of distributions has been proposed to model heavy-tailed count data, but the use of some of them is inhibited by the lack of an explicit, or easily computable, expression for their probability mass function (p.m.f.) and by the lack of any-order moments for all, or some, parameters values.
In this framework, we propose a very general procedure, based on censoring, which requires only the knowledge of the probability generating function (p.g.f.).
One of its appealing characteristics is that parameter estimation is performed by means of a suitably modified moment-based technique, which is appropriate even for distributions without moments. We show that the censored distribution still depends on the parameters of the original one but, having finite moments, allows the application of a moment-based estimation method.
Obviously, the choice of the censoring strength introduces a source of arbitrariness, which can be reduced by adopting data-driven selection criteria. We focus on two-parameter families of distributions and prove that the proposed estimators are consistent and asymptotically normal under rather general mathematical conditions. Moreover, the obtained results can be extended to general multi-parameter families under analogous conditions.

The censoring operation has been already considered for modeling different tail heaviness by generalizing
the Poisson-inverse Gaussian distribution to a more flexible three-parameter family, including, as boundary cases, the Poisson and the discrete stable distributions \citep{Zhu}. The family of discrete stable distributions, introduced by \cite{Steutel}, is large and flexible, allowing skewness, heavy tails, overdispersion, and has many intriguing mathematical properties (see e.g. \citealp{christoph1998discrete}, \citealp{devroye1993triptych}). However, the lack of a closed form expression for the p.m.f. and the non-existence of moments for some parameters values have been a major drawback to its use by practitioners.
Some attempts for parameter estimation have been performed by \cite{Kemp}, \cite{marcheselli2008parameter}, \cite{doray2009some} and \cite{Zhu}, among others. The application of the novel estimation procedure to this family gives rise to estimators  with a closed and simple expression and a really satisfactory performance even for moderate sample size, also when the parameters are in a neighborhood of the boundary values ensuring the existence of moments.  

The general procedure is illustrated in Section \ref{sec:censoring} and results on the discrete stable distributions are given in Section \ref{sec:stable}. Simulation experiments and real data applications are presented in Section \ref{sec:sim-real} while Section \ref{sec:disc} is devoted to concluding remarks. All tables, figures and proofs are reported in the Appendix.

\section{Parameter estimation by censoring}
\label{sec:censoring}
Let $X$ be a count random variable, i.e. a random variable (r.v.) with values in $\mathbb{N}_0=\{0,1,2,\dots\}$ such that $\operatorname{E}[X]$ is not necessarily finite. Moreover, denote by $g$ the p.g.f. of $X$, namely $g(s)=\operatorname{E}[s^X]$ for any $s\in[0,1]$.

When $X$ does not have finite moments and the p.m.f. has no closed-form, any inference procedure becomes cumbersome. To address parameters estimation, we propose an original and effective approach based on  a stochastic perturbation of $X$, giving rise to a censored r.v. with finite moments. More precisely, let $T_p$ be a Geometric r.v. with parameter $p\in ]0,1]$ and p.m.f.
$$h(n)=P(T_p=n)=p(1-p)^{n-1},\,\,\, n\geq 1,$$ 
and let 
$$Y=X\mathbf{1}_{\{X<T_p\}}$$
be the $p$-\emph{censoring} of $X$.
 It is worth noting that $\operatorname{E}[Y]$ is finite since
$$\operatorname{E}[Y]\leq\operatorname{E}[T_p-1]=\frac{1-p}{p}.$$
Moreover, the p.g.f. of $Y$ can be expressed as a function of the p.g.f. of $X$.
\begin{prop}
	\label{pgf_Y}
	Let $g_Y$ be the p.g.f. of $Y$. For any $s\in[0,1]$, 
	\begin{equation}
		\label{gyrelation}
		g_Y(s)=1-g(1-p)+g(s(1-p)).
	\end{equation}
In particular, 
\begin{equation}
\operatorname{E}[Y]=(1-p)g^\prime(1-p),\qquad\operatorname{E}[Y^2]=(1-p)^2g''(1-p)+(1-p)g'(1-p).
\end{equation}
\end{prop}
Now suppose that the distribution of $X$ depends on two parameters $\theta_1,\theta_2$ which can be written as
\begin{equation}\label{teta}
\theta_1=f_1(p,g(1-p),\operatorname{E}[Y]),\quad\theta_2=f_2(p,g(1-p),\theta_1),
\end{equation}
where $f_1$, $f_2$ are two (known) suitable differentiable functions. Condition \eqref{teta} is verified for large classes of distributions, such as discrete stable distributions  {and two-parameter distributions obtained from the Discrete Linnik one fixing the shape parameter}. To estimate $\theta_1$ and $\theta_2$, consider a random sample $(X_1,\dots,X_n)$,  with $X_i\sim X$ for $i=1,\ldots,n$, and $n$ Geometric independent r.v.s $(T_{p,1},\dots,T_{p,n})$ with parameter $p$, independent of $(X_1,\dots,X_n)$,  and let 
\begin{equation}
	\label{empirical}
	\widehat{m}_{p,1}=n^{-1}\sum_{i=1}^{n}X_i\mathbf{1}_{\{  X_i<T_{p,i} \}}, \qquad   \widehat{g}(1-p)=n^{-1}\sum_{i=1}^n(1-p)^{X_i}
\end{equation} 
be the empirical first-order moment of the $p$-censoring of $X$ and the empirical p.g.f. of $X$ computed at $1-p$, respectively. 
The most trivial estimators of $\theta_1$ and $\theta_2$ could be obtained by means of the plug-in technique by replacing $\operatorname{E}[Y]$, that is the finite moment of the $p$-censoring r.v., with its empirical counterpart. Obviously, the variability of the plug-in estimators is inflated by the randomness introduced by censoring. A more precise estimator of $\theta_1$ is provided by considering the conditional expectation of the plug-in estimator given the random sample. Therefore, the proposed estimators are given by
\begin{equation*}\label{stime_con_p}
	\widehat\theta_1=\operatorname{E}[f_1(p,\widehat g(1-p),\widehat{m}_{p,1})|X_1,\ldots,X_n],\quad\widehat\theta_2=f_2(p,\widehat g(1-p),\widehat\theta_1).
\end{equation*}
If a closed form expression for $\widehat\theta_1$ does not exist, it can be approximated by considering $R$ independent generations of $n$ independent Geometric r.v.s with parameter $p$, $(T_{p,1,r},\dots,T_{p,n,r})$, giving rise to $R$ empirical first-order moments $\widehat{m}_{p,1,r}$ ($r=1,\ldots,R)$, in such a way that $\widehat\theta_1$ can be obtained as ${1/R}\sum_{r=1}^Rf_1(p,\widehat g(1-p),\widehat{m}_{p,1,r})$. Obviously, the choice of $R$ is under the control of the researcher and, owing to the negligible computational effort, can be taken large enough to ensure an excellent approximation.
 
Introducing the $p$-censoring r.v. $Y$ induces a source of arbitrariness due to choice of $p$. Indeed, $\widehat\theta_1$ and $\widehat\theta_2$ constitute a family of estimators indexed by $p$ and, thus, the selection of the parameter $p$ ensuring
 enough information on the tail of $X$ and good performance of the corresponding estimators is crucial. In general, values of $p$ in $]0,1/2]$ are advisable. A reasonable approach is to consider a data-driven procedure which should be guided by the features of the p.g.f. of $X$. Then, once the suitable r.v. $p_*$ depending on $X_1,\ldots,X_n$ is defined, the following estimators can be considered
\begin{equation}\label{stime_con_p*}
\widehat\theta_1^*=\operatorname{E}[f_1(p_*,\widehat g(1-p_*),\widehat{m}_{p_*,1})|X_1,\ldots,X_n],\quad\widehat\theta_2^*=f_2(p_*,\widehat g(1-p_*),\widehat\theta_1^*).
\end{equation}
It must be pointed out that, whatever data-driven criterion is adopted, under 
rather mild conditions on $p_*$, thanks to $\eqref{stime_con_p*}$ and to the Delta method, the asymptotic consistency and normality of $\widehat\theta_1^*$ and $\widehat\theta_2^*$ can be proven. 
\begin{prop}
\label{normality}
Suppose there exist $p\in]0,1/2]$ and a sequence $(Z_n)_n$ of independent and identically distributed r.v.s, with $\operatorname{E}[Z_1^2]<\infty$, such that 
\begin{equation}\label{Z}
(p_*-p)-{{\sum_{i=1}^n(Z_i-\operatorname{E}[Z_1])}\over{n}}=o(1)\quad a.s
\end{equation}
and $\sqrt n o(1)=o_P(1)$, where $o_P(1)$ denotes a r.v. which goes to 0 in probability. Moreover, {suppose $f_1$ and $f_2$ be differentiable with respect to $x,y$ and $z$ and denote by ${{\partial f_1}\over{\partial x}}$, ${{\partial f_1}\over{\partial y}}$, ${{\partial f_1}\over{\partial z}}$ and ${{\partial f_2}\over{\partial x}}$, ${{\partial f_2}\over{\partial y}}$, ${{\partial f_2}\over{\partial z}}$ the partial derivatives, respectively. Finally, suppose that} ${{\partial f_1}\over{\partial z}}$ is a bounded function.

\noindent Then $\widehat \theta_1^*$ and $\widehat\theta_2^*$ converge to $\theta_1$ and $\theta_2$ almost surely and 
$[\sqrt n(\widehat \theta_1^*-\theta_1),\sqrt n(\widehat\theta_2^*-\theta_2)]$ converges in distribution to ${\mathcal N}(0,\Sigma)$, where $\Sigma$ is the variance-covariance matrix of $[W_1,W_2]$, with 
$$W_1={{\partial f_1}\over{\partial x}}(P_0)Z_1+{{\partial f_1}\over{\partial y}}(P_0) X^{\prime}_1+{{\partial f_1}\over{\partial z}}(P_0)X^{\prime\prime}_1,$$ 
\begin{multline*}
W_2=({{\partial f_2}\over{\partial x}}(P_1)+{{\partial f_2}\over{\partial z}}(P_1){{\partial f_1}\over{\partial x}}(P_0))Z_1\\
+({{\partial f_2}\over{\partial y}}(P_1)+{{\partial f_2}\over{\partial z}}(P_1){{\partial f_1}\over{\partial y}}(P_0)) X^{\prime}_1+{{\partial f_2}\over{\partial z}}(P_1){{\partial f_1}\over{\partial z}}(P_0)X^{\prime\prime}_1
\end{multline*}
and 
\begin{equation}\label{X_primo}
	 X^{\prime}_1=(1-p)^{X_1}- \operatorname{E}[X_1(1-p)^{X_1-1}]Z_1, \ \ X^{\prime\prime}_1=X_1(1-p)^{X_1}-\operatorname{E}[X_1^2(1-p)^{X_1-1}]Z_1,  
\end{equation}
$$P_0=(p,g({1-}p),\operatorname{E}[X_1(1-p)^{X_1}]), \ \ P_1=(p,g({1-}p),\theta_1).$$ 
\end{prop}
{Condition \eqref{Z} requires that the sequence $(p_*-p)_n$  is asymptotically equivalent to a sequence of averages of i.i.d. centered r.v.s. Any data driven criterion, giving rise to a $p^*$ which satisfies \eqref{Z}, ensures the asymptotic properties of $[\widehat \theta_1, \widehat \theta_2]$.} Moreover, when the parameter of the Geometric r.v. is not selected by a data-driven procedure but it is fixed in advance, the asymptotic properties of $[\widehat \theta_1, \widehat \theta_2]$ hold under the sole assumption that ${{\partial f_1}\over{\partial z}}$ is bounded.

As to the estimation of $\Sigma$, a suitable estimator  is given by the sample variance-covariance matrix of $[W_{1,1}^*,W_{2,1}^*],\ldots,[W_{1,n}^*,W_{2,n}^*]$, where
$$W_{1,i}^*={{\partial f_1}\over{\partial x}}(P^*_0)Z_i+{{\partial f_1}\over{\partial y}}(P^*_0) X^{\prime}_i+{{\partial f_1}\over{\partial z}}(P^*_0)X^{\prime\prime}_i,$$  
\begin{multline*}
	W_{2,i}^*=({{\partial f_2}\over{\partial x}}(P_1^*)+{{\partial f_2}\over{\partial z}}(P_1^*){{\partial f_1}\over{\partial x}}(P_0^*))Z_i\\
	+({{\partial f_2}\over{\partial y}}(P^*_1)+{{\partial f_2}\over{\partial z}}(P^*_1){{\partial f_1}\over{\partial y}}(P^*_0)) X^{\prime}_i+{{\partial f_2}\over{\partial z}}(P^*_1){{\partial f_1}\over{\partial z}}(P^*_0)X^{\prime\prime}_i,
\end{multline*}
$$P^*_0=(p_*,\widehat g({1-}p_*),\widehat m_{p_*,1}), \quad P^*_1=(p_*,\widehat g({1-}p_*),\widehat\theta_1^*),$$
 and $ X^{\prime}_i$, $X^{\prime\prime}_i$ have expressions analogous to those in \eqref{X_primo}. Indeed, under the assumptions of the previous proposition, $P^*_0$, $P^*_1$ converge almost surely to $P_0,P_1$ and therefore the sample variance-covariance matrix is a consistent estimator of $\Sigma$.

\section{Parameter estimation for the discrete stable family}
\label{sec:stable}
The discrete stable family, denoted as $\mathcal{DS}(a,\lambda)$ with $a\in ]0,1]$ and $\lambda>0$, constitutes an interesting two-parameter model on $\mathbb{N}_0$ with a Paretian tail, whose use is inhibited by the lack of an explicit expression for its p.m.f. and of moments of any order when $a<1$. Indeed, these features preclude the exploitation of the maximum-likelihood or the moment method for parameter estimation. However, since its p.g.f. is 
\begin{equation}
	g(s)=\exp\left( -\lambda(1-s)^a\right),
\end{equation}
the proposed censoring technique can be suitably applied.
Obviously, for $a=1$ the discrete stable family reduces to the Poisson family of distributions and $\operatorname{E}[X]=\infty$ when $a<1$.
Now, let $X_1,\dots,X_n$ be a random sample from $X\sim\mathcal{DS}(a,\lambda)$. 
By using \eqref{gyrelation}, the p.g.f of the $p$-censoring turns out to be
$$g_Y(s)=1-e^{-\lambda p^a}+e^{-\lambda(1-s(1-p))^a},$$
in such a way that
\begin{equation}\label{per_lambda}
	\operatorname {E}[Y]=g(1-p)p^{a-1}(1-p)\lambda a
\end{equation}
and, by noting that $g(1-p)=\exp\left( -\lambda p^a\right)$,  
\begin{equation}\label{per_a}
\operatorname {E}[Y]=-ap^{-1}(1-p)g(1-p)\log(g(1-p)).
\end{equation}
From \eqref{per_a} and \eqref{per_lambda}, it is at once apparent that $a$ and $\lambda$ can be expressed as in \eqref{stime_con_p} and thus, from \eqref{stime_con_p*}, they can be estimated by means of 
\begin{align}
	\label{estimator_a}
	\widehat a \nonumber &=-\operatorname{E}[{{p_*\widehat m_{p_*,1}}\over {(1-p_*)\widehat g(1-p_*)\, \log\big( \widehat g(1-p_*) \big)}}|X_1,\ldots,X_n]\\ \nonumber
	&=-{{p_*\operatorname{E}[m_{p_*,1}|X_1,\ldots,X_n]}\over {(1-p_*)\widehat g(1-p_*)\, \log\big( \widehat g(1-p_*) \big)}}\\
	&=-{{p_*}\over {(1-p_*)\widehat g(1-p_*)\, \log\big( \widehat g(1-p_*) \big)}}n^{-1}\sum_{i=1}^n X_i(1-p_*)^{X_i}
\end{align}
 \begin{equation}
		\label{estimator_lambda}
\widehat \lambda=-p_*^{-\widehat a}\log\big(\widehat g(1-p_*) \big),
\end{equation}
where $p_*$, representing the data-driven choice of the censoring parameter, depends only on $X_1,\ldots, X_n$.

In order to ensure satisfactory finite sample performance of estimators $\eqref{estimator_a}$ and \eqref{estimator_lambda}, the data-driven choice of $p_*$ is crucial. In particular, $p_*$ should be chosen to provide that the denominator in \eqref{estimator_a} is not too close to zero. Then, since $p_*$ is less or equal to $1/2$ and $x\mapsto-x\log x$ has maximum at $x=1/e$, we propose the following data-driven criterion
\begin{equation}
	\label{p_star}
p_*={\rm argmax}_{p\in]0,1/2]}\{-\widehat g(1-p)\log\big(\widehat g(1-p) \big)\}=\max\{p\in]0,1/2]:\widehat g(1-p)\geq 1/e\}.
\end{equation}
Note that if $\widehat g(1/2)\geq 1/e$, $p_*=1/2$ and \eqref{estimator_a} reduces to
$$	\widehat a = -{{n^{-1}\sum_{i=1}^n X_i2^{-X_i}}\over {\widehat g(1/2)\, \log\big( \widehat g(1/2) \big)}}$$
while if $\widehat g(1/2)< 1/e$, $p_*$ is less than $1/2$ and such that $\widehat g(1-p)=1/e$ and the estimator of $a$ is given by
$$	\widehat a = {{ep_*n^{-1}\sum_{i=1}^n X_i(1-p_*)^{X_i}}\over {1-p_*}}.$$
Then, the estimators for $a$ and $\lambda$ turn out to be
\begin{equation}
	\label{a_pstar}
	\widehat a ={{ep_*n^{-1}\sum_{i=1}^n X_i(1-p_*)^{X_i}}\over {1-p_*}}\mathbf{1}_{\{p_*<1/2\}}-{{n^{-1}\sum_{i=1}^n X_i2^{-X_i}}\over {\widehat g(1/2)\, \log\big( \widehat g(1/2) \big)}}\mathbf{1}_{\{p_*=1/2\}}
\end{equation} 
\begin{equation}
	\label{lambda_pstar}
	\widehat \lambda=p_*^{-\widehat a}\mathbf{1}_{\{p_*<1/2\}}-2^{\, \widehat a}\log\big( \widehat g(1/2) \big)\mathbf{1}_{\{p_*=1/2\}}.
\end{equation}
It is worth noting that asymptotically $p_*<1/2$  almost surely when $\lambda >2^a$, otherwise $p_*=1/2$ if $\lambda<2^a$. 
\begin{prop}
	\label{pconvergence}
	$p_*$ converges almost surely to $p=\min(\lambda^{-1/a},1/2)$. Moreover, for $\lambda>2^a$ condition $(6)$ holds with $Z_i=e\lambda^{-1/a}(1-\lambda^{-1/a})^{X_i}/a$ while, for $\lambda<2^a$, with  $Z_1=0$. 
\end{prop}
By using Proposition \ref{normality} and Proposition \ref{pconvergence}, consistency and asymptotic normality of estimators $\widehat{a},\, \widehat{\lambda}$ are obtained.
\begin{cor}
	\label{cor}
	{\it $\widehat{a},$ $\widehat{\lambda}$ converge almost surely to $a,\lambda$ and $[\sqrt n(\widehat a-a),\sqrt n(\widehat\lambda-\lambda)]$ converges in distribution to ${\mathcal N}(0,\Sigma)$, where $\Sigma$ is the variance-covariance matrix of $[W_1,W_2]$, with  $$W_1=e\lambda^{-1/a}
		X_1(1-\lambda^{-1/a})^{X_1-1},$$ $$ W_2=-e\lambda\big((1-\lambda^{-1/a})^{X_1}-a^{-1}(\lambda^{-1/a}
		\log {\lambda}) X_1(1-\lambda^{-1/a})^{X_1-1}\big),$$ 
		for $\lambda>2^a $, while, for $\lambda<2^a$, 
		$$W_1=\lambda^{-1}2^a
		e^{{\lambda}\over{2^a}}(X_12^{-X_1}+a(1-\lambda 2^{-a})2^{-X_1}),$$
		$$ W_2=2^a
		e^{{\lambda}\over{2^a}}(X_12^{-X_1}\log 2+\big(a(1-\lambda 2^{-a})\log 2-1\big)2^{-X_1}).$$}
\end{cor}
As in the general case  in Section \ref{sec:censoring}, a suitable estimator of $\Sigma$ is given by the sample variance-covariance matrix of $[W_{1,1}^*,W_{2,1}^*],\ldots,[W_{1,n}^*,W_{2,n}^*]$, where
 $$W_{1,i}^*=ep_*
 X_i(1-p_*)^{X_i-1},$$
 $$ W_{2,i}^*=-e{\widehat\lambda}\big((1-p_*)^{X_i}+{X_i(1-p_*)^{X_i-1}}p_*
 \log{p_*}\big),$$ for  $p_*<1/2$ while, for $p_*=1/2$, 
$$W_{1,i}^*=-{{2^{-X_1}
\big(X_1+\widehat a(1+\log\widehat g(1/2))\big)}\over{\widehat g(1/2)\log\widehat g(1/2)}},$$ $$ W_2=2^{\widehat{a}-X_1}
 e^{{\widehat\lambda}\over{2^{\widehat a}}}\big(X_1\log 2+\big(\widehat a(1-\widehat\lambda 2^{-\widehat a})\log 2-1\big)\big).$$

\section{Simulation study and real data applications}
\label{sec:sim-real}
\subsection{Simulation study}
The performance of the estimators \eqref{a_pstar} and  \eqref{lambda_pstar} was assessed by means of an extensive
Monte Carlo simulation implemented by using R (\citealp{R2020}). Following \cite{devroye1993triptych}, the realizations of the discrete stable distribution were generated by using the equality in law $\mathcal{DS}(a, \lambda)\overset{\mathcal{L}}{=} \mathcal{P}(\mathcal{PS}(a,\lambda)),$
where {$\mathcal{P}(\mathcal{PS}(a,\lambda))$ denotes a Poisson compound probability distribution where the Poisson parameter is a random variable with positive stable distribution}.
For generating realizations from the positive stable distribution, the classical Kanter’s representation (\citealp{Kanter})
was adopted. 
As to the parameter values, the value of $a$ was set equal to $0.25, 0.5, 0.75, 1$ while all the values varying from $0.5$ to $12$ by $0.5$ were considered for $\lambda$. 

\begin{figure} [H]
	\centering
	\includegraphics[width=\textwidth]{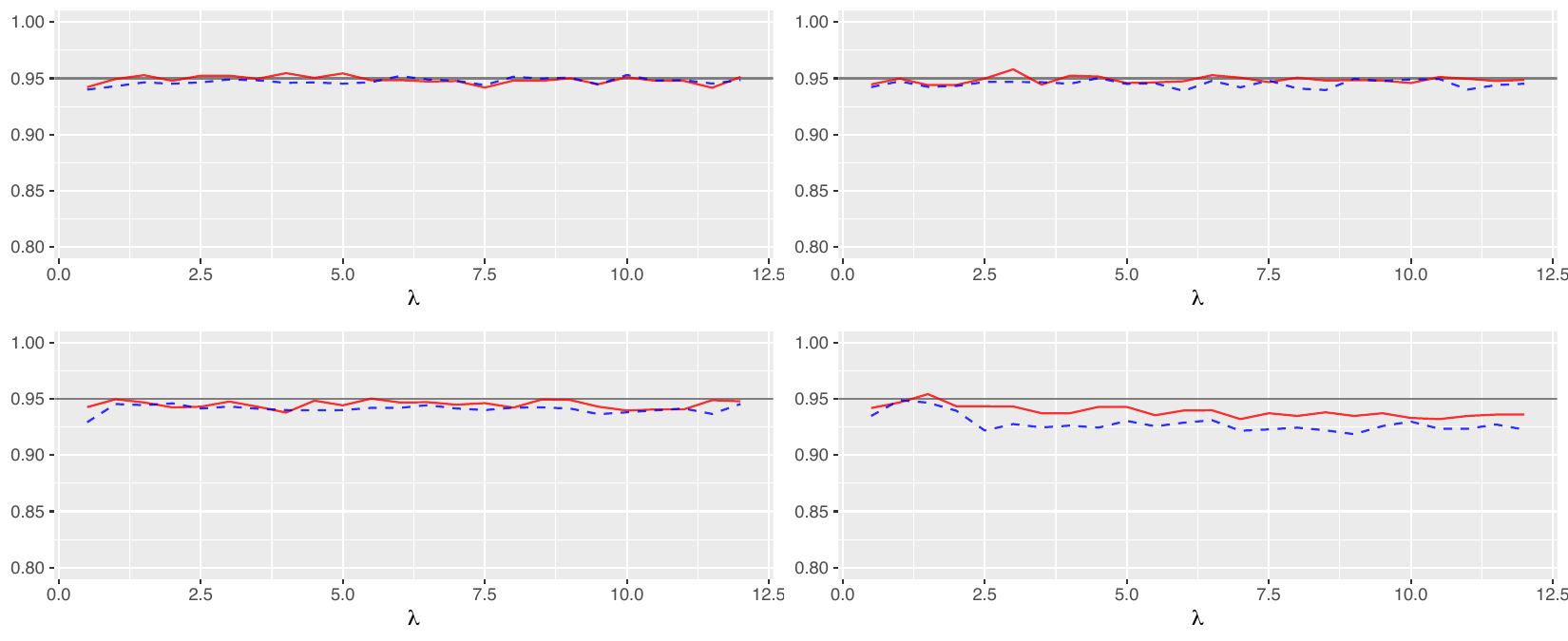}
	\caption{Empirical coverage of 95\% confidence intervals for ${a}$ with $n=100$ (dashed line) and $n=200$ (solid line). Top-left: $a=0.25$, top-right: $a=0.5$, bottom-left: $a=0.75$, bottom-right: $a=1$.}
	\label{fig:figure1}
\end{figure}

\begin{figure}[H] 
	\centering
	\includegraphics[width=\textwidth]{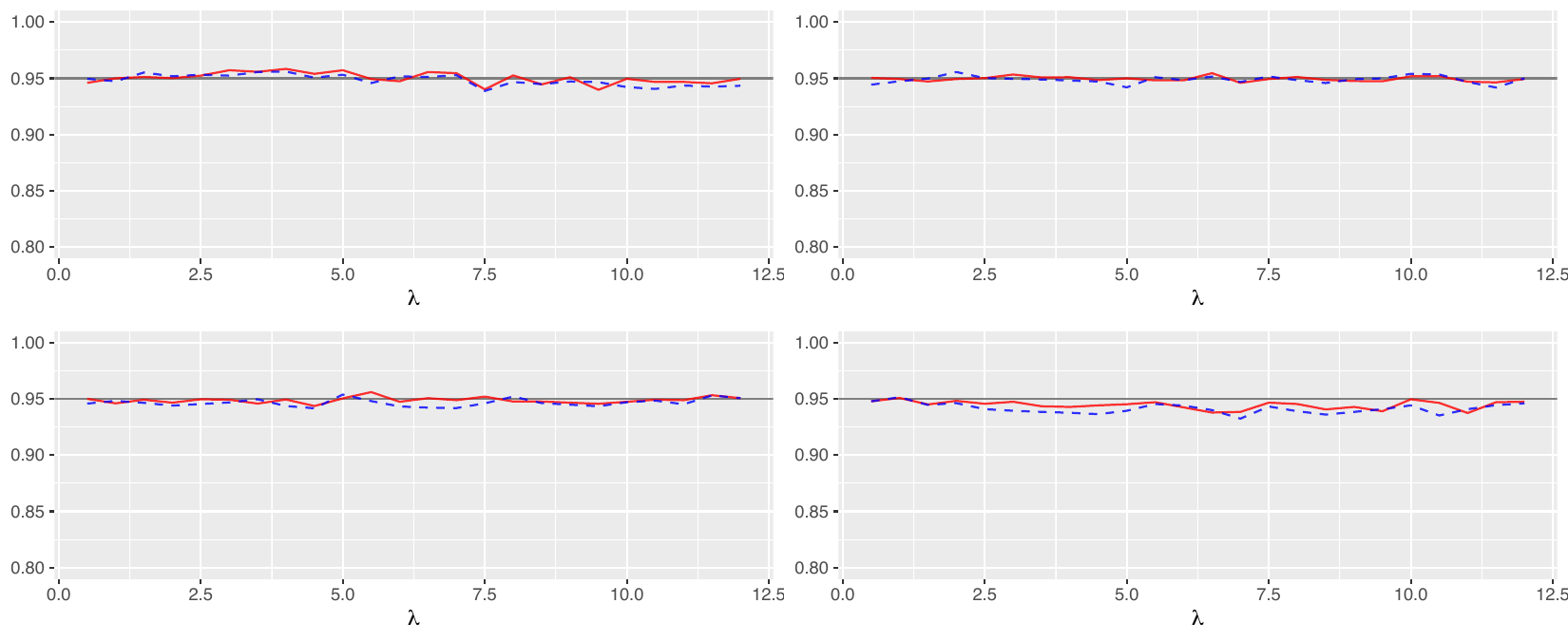}
	\caption{Empirical coverage of 95\% confidence intervals for ${\lambda}$ with $n=100$ (dashed line) and $n=200$ (solid line). Top-left: $a=0.25$, top-right: $a=0.5$, bottom-left: $a=0.75$, bottom-right: $a=1$.}
	\label{fig:figure2}
\end{figure}

For each combination of $a$ and $\lambda$ values, 
$5000$ samples of size $n=100,200$ were independently generated from $\mathcal{DS}(a, \lambda)$. For each sample, first the censoring parameter was selected according to \eqref{p_star} and then parameter estimates were obtained by means of \eqref{a_pstar} and  \eqref{lambda_pstar}, together with the corresponding variance estimates. Moreover, confidence interval estimates for $a$ and $\lambda$ at confidence level $0.95$ were obtained using the quantiles of the standard normal distribution.  From the Monte Carlo distributions, the Relative Root Mean Squared Error (RRMSE) of \eqref{a_pstar} and  \eqref{lambda_pstar} was obtained and reported in Table \ref{tab:tabella1} and Table \ref{tab:tabella2}, respectively. For any combination of $\lambda$ and $a$, the RRMSEs of both estimators are rather satisfactory and obviously decrease as $n$ increases. Moreover, for any fixed $n$ and $a$, the RRMSE of $\widehat a$ decreases as $\lambda$ increases. Similarly, for any fixed $n$ and $\lambda$, the RRMSE of $\widehat \lambda$ decreases as $a$ increases for the larger values of $\lambda$. {T}he empirical coverages of the confidence intervals were {also} computed. For any fixed value $a=0.25, 0.5, 0.75, 1$, Figures \ref{fig:figure1} and \ref{fig:figure2} depict the empirical coverage of the 0.95 confidence intervals for $a$ and $\lambda$ respectively, with $\lambda$ varying from $0.5$ to $12$ by $0.5$, and both for $n=100$ and $n=200$.
From Figures \ref{fig:figure1} and \ref{fig:figure2} it is apparent that the empirical coverages are really satisfactory for both parameters even for $n=100$. Empirical coverages of the confidence intervals for $a$ are less close to $0.95$ only when $a=1$ and $n=100$, while for $n=200$ they approach the nominal one. {Finally, we also performed simulations for fixed $p$ whose results, not reported for the sake of brevity, strongly support the proposed data-driven procedure.}

\subsection{Real data application}
We fit the discrete stable distribution on citation data from Web of Science database. In particular, the dataset was composed of 369 citation counts of articles published in 2000 with keyword \lq\lq linear model". We obtained $\widehat{a}= 0.5583$, $\widehat{\lambda}= 4.4689$
and $95\%$ interval estimates $[0.5188,0.5977]$ and $[3.894,5.0435]$ at confidence level 0.95 for $a$ and $\lambda$, respectively. 

\begin{figure}[H]
	\centering
	\includegraphics[width=0.45\textwidth]{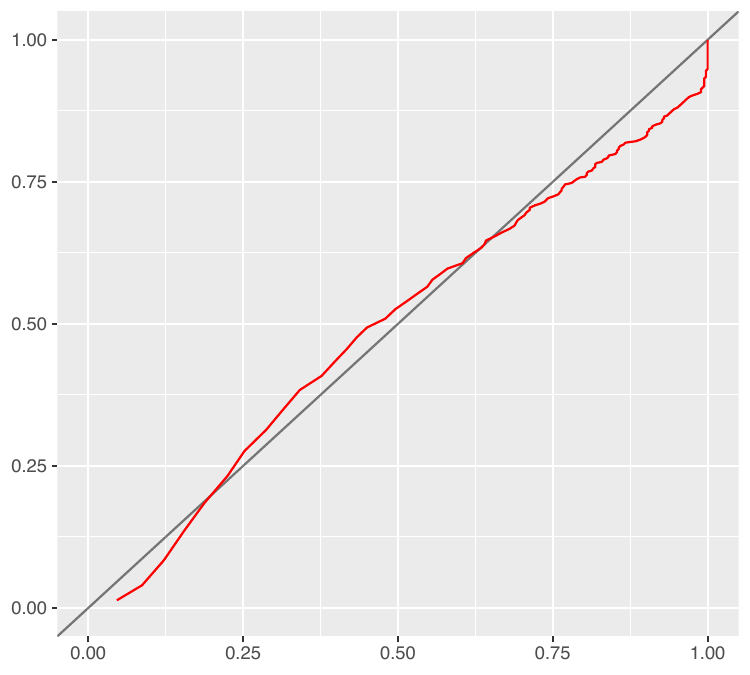}
	\caption{$F_n \circ Q_X$ ($x$-axis) and $F_X \circ Q_X$ ($y$-axis), where $F_X$ and $Q_X$ are respectively the theoretical distribution function and theoretical quantile function (computed by simulation), while $F_n$ is the empirical distribution function.  }
	\label{fig:wos}	
\end{figure}

\begin{figure}[H]
	\centering
	\includegraphics[width=0.45\textwidth]{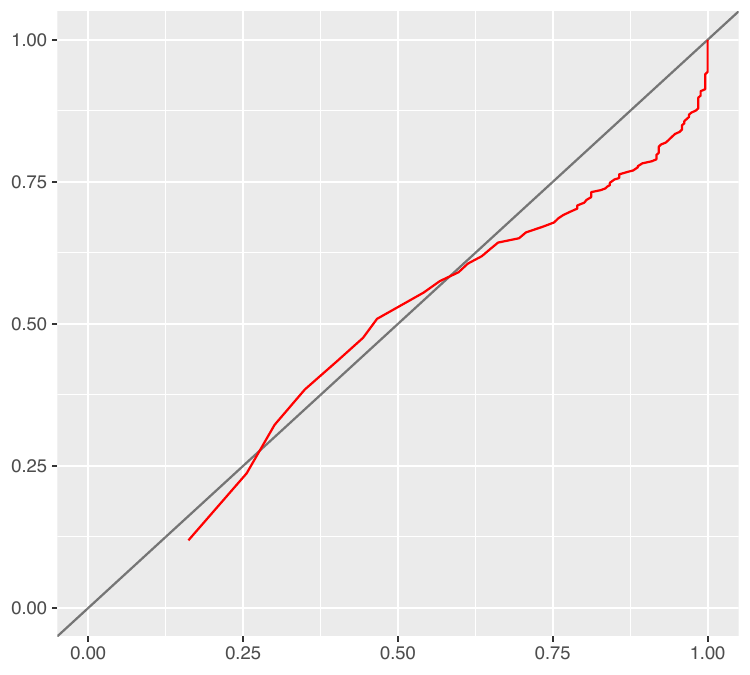}
	\caption{$F_n \circ Q_X$ ($x$-axis) and $F_X \circ Q_X$ ($y$-axis), where $F_X$ and $Q_X$ are respectively the theoretical distribution function and theoretical quantile function (computed by simulation), while $F_n$ is the empirical distribution function.}
	\label{fig:Zhou}
\end{figure}

We also considered the citation data presented in \cite{Zhu} and obtain $\widehat{a}= 0.4613$, $\widehat{\lambda}= 2.1471$
and $95\%$ interval estimates $[0.4104, 0.5123]$ and $[1.8203,2.4738]$ for $a$ and $\lambda$, respectively. 
Graphical assessment of the model's goodness-of-fit for both datasets is reported in Figures \ref{fig:wos} and \ref{fig:Zhou}.

\section{Discussion}
\label{sec:disc}
A general procedure for parameter estimation is welcomed when the maximum likelihood or moment-based criteria are precluded and, even more so, if it allows to obtain consistent and asymptotically normal estimators. The proposed procedure, under mild mathematical conditions, not only gives rise to estimators sharing these properties, but also to variance estimators which avoid computationally intensive resampling methods. For the discrete stable family, the proposed estimators also show rather satisfactory performance for finite sample, in terms of coverages of confidence intervals and of relative root mean squared errors. The novel estimation procedure has been introduced referring to distributions depending on two parameters, but it could be generalized to distributions with more parameters. Obviously, when $k$ parameters are under estimation, generally moments up to the order $k-1$ for the $p$-censoring r.v. are involved and they can be straightforwardly obtained by means of Proposition \ref{pgf_Y}.
Finally, further research will be devoted to investigate if the censoring could also be used to introduce a general goodness-of-fit test for count distributions (also without any moments) as few proposals are available and many of them, being tailored to deal with particular distributions, are of limited applicability.

\appendix
\section*{Appendix}

\section{Tables}

\begin{table}[h]
	 \centering
	\caption{Percentage values of the RRMSE of $\widehat a$ for various combinations of $\lambda$ values, $a$ values and sample sizes.}
    \medskip
\begin{tabular}{lcccccccc} 
		\toprule
			& \multicolumn{2}{c}{$a=0.25$}& \multicolumn{2}{c}{$a=0.5$}&\multicolumn{2}{c}{$a=0.75$}&\multicolumn{2}{c}{$a=1$}\\	
		 	 \cmidrule(lr){2-3} \cmidrule(lr){4-5}  \cmidrule(lr){6-7} \cmidrule(lr){8-9}
			$\lambda$ &$n=100$& $n=200$& $n=100$& $n=200$ & $n=100$ & $n=200$ & $n={100}$ &$n={200}$\\
	    \midrule
			$0.5$        &  23 &16 &15 &10& 10& 7& 4& 3          \\
			$1$             & 20& 14 &12& 9& 8& 6& 4& 3          \\
			$2$             & 14 &10& 9& 7& 7& 5& 5& 3           \\
			$5$          &  13& 10& 8& 5& 5& 3& 2& 1               \\
			$10$           & 14& 9& 7& 5& 4& 3& 1& 1             \\
			
			\bottomrule
		\end{tabular}  
	\label{tab:tabella1}
\end{table}

\begin{table}[h]
	 \centering
	\caption{Percentage values of the RRMSE of $\widehat\lambda$ for various combinations of $\lambda$ values, $a$ values and sample sizes.}
	 \medskip
 \begin{tabular}{lcccccccc} 
         \toprule
			& \multicolumn{2}{c}{$a=0.25$}& \multicolumn{2}{c}{$a=0.5$}&\multicolumn{2}{c}{$a=0.75$}&\multicolumn{2}{c}{$a=1$}\\	
	 		 \cmidrule(lr){2-3} \cmidrule(lr){4-5}  \cmidrule(lr){6-7} \cmidrule(lr){8-9}			
			$\lambda$ &$n=100$& $n=200$& $n=100$& $n=200$ & $n=100$ & $n=200$ & $n={100}$ &$n={200}$\\
   \midrule		
			$0.5$            & 	16 & 11 & 16 & 11 & 16  &11 &16 & 11            \\
			$1$               & 13&9& 13& 9   &13   &9   &12  & 8         \\
			$2$                &   14 &10& 11& 8& 10& 7& 10& 7             \\
			$5$             &    25& 16 &13 &9 &8 &6 &6 &4                \\
			$10$             &       36& 23& 17& 11& 9 &6 &4& 3                 \\
		\bottomrule
		\end{tabular}  
	\label{tab:tabella2}
\end{table}

\section{Figures}

\section{Proofs}

\subsection*{Proof of Proposition \ref{pgf_Y}}
Since $P(T_p\leq n)=1-(1-p)^n$, for $n\geq 1$ it holds 
$$P(Y=n)=P(X=n)P(T_p> n)=P(X=n)(1-p)^n.$$
Moreover, it must be pointed out that 
\begin{equation}
	\label{TgraeterX}
	g(1-p)=\sum_{n=0}^{\infty}(1-p)^nP(X=n)=P(T_p>X).
\end{equation}
Thus
\begin{align*}
	g_Y(s)&=\sum_{n=0}^\infty s^nP(Y=n)\\
	&=P(Y=0)+\sum_{n=1}^\infty s^n(1-p)^nP(X=n)\\
	&=P(Y=0)-P(X=0)+g\big(s(1-p)\big)\\
	&=P(T_p\leq X)+g\big(s(1-p)\big)\\
	 &=1-g(1-p)+g\big(s(1-p)\big).
\end{align*}
\noindent
From the previous expression, it is at once apparent that $$\operatorname{E}[Y]=g^\prime_Y(1)=(1-p)g^\prime(1-p)$$ and 
$$\operatorname{E}[Y^2]=g^{\prime\prime}_Y(1)+\operatorname{E}[Y]=(1-p)^2g''(1-p)+(1-p)g'(1-p),$$
and the proof is concluded.\qed

\subsection*{Proof of Proposition \ref{normality}} 
Let $L$ be a constant such that $|{{\partial f_1}\over{\partial z}} |\leq L$. Since
\begin{align*}
	\Big|\widehat\theta_1^*&-f_1(p_*,\widehat g(1-p_*),\frac{1}{n}\sum_{i=1}^nX_i(1-p_*)^{X_i})\Big|\\
	&\leq L\operatorname{E}[(\widehat{m}_{p_*,1}-\frac{1}{n}\sum_{i=1}^nX_i(1-p_*)^{X_i})^2|X_1,\ldots,X_n ]^{1\over 2}
\end{align*}
and
\begin{align*}
	\operatorname{E}[(\widehat{m}_{p_*,1}-\frac{1}{n}\sum_{i=1}^nX_i(1-p_*)^{X_i})^2|X_1,\ldots,X_n ]&=
	\frac{1}{n^2}\sum_{i=1}^n X_i^2(1-p_*)^{X_i}\big(1-(1-p_*)^{X_i}\big)\\
	&=O(1/n) \quad a.s.,
\end{align*}
 the consistency of $\widehat \theta_1^*$ and, consequently, of $\widehat\theta_2^*$ is obtained if $$\lim_n 
f_1(p_*,\widehat g(1-p_*),\frac{1}{n}\sum_{i=1}^nX_i(1-p_*)^{X_i})=\theta_1\quad a.s.$$ holds. In other words, thanks to continuity of $f_1$ and to condition \eqref{Z}, which implies $\lim_n p_*=p$ a.s., it suffices to prove $$\lim_n\widehat g(1-p_*)=g(1-p),\qquad \lim_n\frac{1}{n}\sum_{i=1}^nX_i(1-p_*)^{X_i}=\operatorname{E}[X_1(1-p)^{X_1}]\qquad a.s.$$ From the Strong Law of Large Numbers and from 
\begin{equation}\label{Taylor}
	|(1-p_*)^{X_i}-(1-p)^{X_i}|\leq X_i2^{1-X_i}|p_*-p|,
\end{equation}
the previous relations, and then the consistency, immediately follow. 

Moreover, from \eqref{Z} and \eqref{Taylor}, by applying again the Strong Law of Large Numbers, it holds
\begin{align*}
	\sqrt n\Big(&{1\over n}\sum_{i=1}^n X_i(1-p_*)^{X_i}-{1\over n}\sum_{i=1}^n X_i(1-p)^{X_i}\Big)\\
	&=-\operatorname{E}[X_1^2(1-p)^{X_1-1}]\sqrt n(p_*-p)+o(1)\\
	&=-\operatorname{E}[X_1^2(1-p)^{X_1-1}]{{\sum_{i=1}^n(Z_i-\operatorname{\operatorname{E}}[Z_1])}\over{\sqrt n}}+o(1),
\end{align*}
which implies
\begin{equation} \label{Xprimeprime}
	\sqrt n\Big({1\over n}\sum_{i=1}^n X_i(1-p_*)^{X_i}-\operatorname{E}[X_1(1-p)^{X}]\Big)={{\sum_{i=1}^n(X^{\prime\prime}_i-\operatorname{E}[X_1^{\prime\prime}])}\over{\sqrt n}}+o(1),
\end{equation} 
where $X^{\prime\prime}_i=X_i(1-p)^{X_i}- \operatorname{E}[X_1^2(1-p)^{X_1-1}]Z_i$. Similarly,
\begin{equation}\label{Xprime}
	\sqrt n\big(\widehat g(1-p_*)-g(1-p)\big)={{\sum_{i=1}^n (X^{\prime}_i-\operatorname{E}[X^{\prime}_1])}\over{\sqrt n}}+o(1),
\end{equation}
where $X^{\prime}_i=(1-p)^{X_i}- \operatorname{E}[X_1(1-p)^{X_1-1}]Z_i$.

Now, let $P_0=(p,g(p),\operatorname{E}[X_1(1-p)^{X_1}])$. From \eqref{Z}, \eqref{Xprimeprime} and \eqref{Xprime}, it follows 
\begin{align}\label{distthetahat1}
	\begin{split}
    &\sqrt n(\widehat \theta_1^*-\theta_1)\\
	&={{\sum_{i=1}^n {{\partial f_1}\over{\partial x}}(P_0)(Z_i-\operatorname{E}[Z_1])+{{\partial f_1}\over{\partial y}}(P_0) (X^{\prime}_i-\operatorname{E}[X^{\prime}_1])+{{\partial f_1}\over{\partial z}}(P_0)(X^{\prime\prime}_i-\operatorname{E}[X^{\prime\prime}_1])}\over{\sqrt n}}+o_P(1).
\end{split}
\end{align} 
Owing to the classical Central Limit Theorem, $\sqrt n(\widehat \theta_1^*-\theta_1)$ converges in distribution to ${\mathcal N}(0,\sigma_1^2)$ where
\begin{equation*}
	\sigma_1^2={\rm Var}[{{\partial f_1}\over{\partial x}}(P_0)Z_1+{{\partial f_1}\over{\partial y}}(P_0) X^{\prime}_1+{{\partial f_1}\over{\partial z}}(P_0)X^{\prime\prime}_1].
\end{equation*} 
Finally, by arguing in a similar way, from \eqref{distthetahat1} it follows
\begin{align}\label{distthetahat2}
	\begin{split}
&\sqrt n(\widehat \theta_2^*-\theta_2)\\
	&={{\sum_{i=1}^n {{\partial f_2}\over{\partial x}}(P_1)(Z_i-\operatorname{E}[Z_1])+{{\partial f_2}\over{\partial y}}(P_1) (X^{\prime}_i-\operatorname{E}[X^{\prime}_1])}\over{\sqrt n}}+{{\partial f_2}\over{\partial z}}(P_1)\sqrt n(\widehat \theta_1-\theta_1)+o_P(1)
\end{split}
\end{align} 
where $P_1=(p,g(1-p),\theta_1)$. In particular, $\sqrt n(\widehat \theta_2^*-\theta_2)$ converges in distribution to ${\mathcal N}(0,\sigma_2^2)$ where
\begin{multline*}
	\sigma_2^2={\rm Var}[({{\partial f_2}\over{\partial x}}(P_1)+{{\partial f_2}\over{\partial z}}(P_1){{\partial f_1}\over{\partial x}}(P_0))Z_1\\
	+({{\partial f_2}\over{\partial y}}(P_1)+{{\partial f_2}\over{\partial z}}(P_1){{\partial f_1}\over{\partial y}}(P_0)) X^{\prime}_1+{{\partial f_2}\over{\partial z}}(P_1){{\partial f_1}\over{\partial z}}(P_0)X^{\prime\prime}_1].
\end{multline*} 
Thus, from \eqref{distthetahat1} and \eqref{distthetahat2} the thesis follows.\qed
\subsection*{Proof of Proposition \ref{pconvergence}}
Let $F_n(x)={{1}\over{n}}\sum_{i=1}^n I_{\{X_i\leq x\}}$ and $F$ be the distribution function of $X$. It holds
$$\sup_{p\in ]0,1/2]}|\widehat g(1-p)-g(1-p)|\leq \sup_x \big|\widehat F_n(x)-F(x)\big|$$ 
(see e.g., \citealp{marcheselli2008parameter}, page 824).
Therefore, owing to Glivenko-Cantelli Theorem it follows
 $$\sup_{p\in ]0,1/2]}|\widehat g(1-p)-g(1-p)| \to 0, \ \ a.s. $$ In particular, if $\lambda^{-1/a}>1/2$, then asymptotically it holds $p^*=1/2$ a.s. Moreover, if $\lambda^{-1/a}<1/2$, then  $\widehat{g}(1-p_*)=e^{-1}$ a.s.\ for large $n$ and  
\begin{align*}\label{g}
	g(1-{\lambda}^{-{1\over a}})-\widehat g(1-{\lambda}^{-{1\over a}})
	&=1/e-\widehat g(1-{\lambda}^{-{1\over a}})\\
	&=\widehat{g}(1-p_*)-\widehat g(1-{\lambda}^{-{1\over a}}).
\end{align*}
Thanks to Lagrange Theorem, for large $n$, there exists $C_{n}\in ]\min(p_*,{\lambda}^{-{1\over a}}),\max(p_*,{\lambda}^{-{1\over a}})[$ such that
\begin{equation}\label{g}
	g(1-{\lambda}^{-{1\over a}})-\widehat g(1-{\lambda}^{-{1\over a}})=-\widehat g^\prime(1-C_{n})(p_*-{\lambda}^{-{1\over a}})\ \ \ \ \ \  a.s.
\end{equation}
Since $\widehat g^\prime$ is non-decreasing it follows
\begin{align}\label{gprime}
	\begin{split}
	\widehat g^\prime(1/2)&\leq\widehat g^\prime(1-\max(p_*,{\lambda}^{-{1\over a}}))\\
	&\leq \widehat g^\prime(1-C_{n})\\
	&\leq \widehat g^\prime(1-\min(p_*,{\lambda}^{-{1\over a}}))
		\end{split}
\end{align}
and $$\lim_n \widehat g(1-{\lambda}^{-{1\over a}})=g(1-{\lambda}^{-{1\over a}})\quad a.s.$$
Thus, from \eqref{g} and \eqref{gprime}, $p_*$ converges a.s.\ to $ {\lambda}^{-{1\over a}}$. (Note that $p_*$ converges a.s.\ to ${1/2}$ when ${\lambda}^{-{1\over a}}=1/2$). Moreover, from \eqref{gprime} it holds

\begin{equation}\label{gprime2}
	\widehat g^\prime(1-C_{n}) \overset{a.s.}{\to}g^\prime(1-\lambda^{-{1\over a}})={{a\lambda^{1\over a}}\over{e}}. 
\end{equation}
From \eqref{g} and \eqref{gprime2}
\begin{align*}
	(p_*-{\lambda}^{-{1\over a}})&=\frac{g(1-{\lambda}^{-{1\over a}})-\widehat g(1-{\lambda}^{-{1\over a}})}{-\widehat g^\prime(1-C_{n})}\\
	&\sim{{e{\lambda}^{-{1\over a}}}\over a}\frac{\sum_{i=1}^n(1-{\lambda}^{-{1\over a}})^{X_i}-\operatorname{E}[(1-{\lambda}^{-{1\over a}})^{X_1}]}{n} \ \ a.s.
\end{align*}
and the thesis follows.\qed
\subsection*{Proof of Corollary \ref{cor}}
Let $\lambda>2^a$ and $p={\lambda}^{-{1\over a}}$. In this case, $f_1,f_2$ are defined by
$$f_1(x,y,z)={{exz}\over{1-x}},\qquad f_2(x,y,z)=x^{-z},$$
$|{{\partial f_1}\over{\partial z}}|\leq e$ and condition $\eqref{Z}$ of Proposition \ref{normality} is verified with $$Z_i=\frac{ep(1-p)^{X_i}}{a}.$$ Then, $\widehat a$ and $\widehat\lambda$ converge a.s.\ to $a$ and $\lambda$. Moreover, since $$\operatorname{E}[X_1(1-p)^{X_1}]=p\operatorname{E}[X_1^2(1-p)^{X_1}],$$ from Proposition \ref{normality}, after a little algebra, $[\sqrt n(\widehat a-a),\sqrt n(\widehat\lambda-\lambda)]$ converges in distribution to ${\mathcal N}(0,\Sigma)$, where $\Sigma$ is the variance-covariance matrix of $[W_1,W_2]$, with  $$W_1={{e{\lambda}^{-{1\over a}}}\over{1-{\lambda}^{-{1\over a}}}}X_1(1-{\lambda}^{-{1\over a}})^{X_1}, $$ 
$$ W_2=-e{\lambda}\big((1-{\lambda}^{-{1\over a}})^{X_1}+({\lambda}^{-{1\over a}}\log {\lambda}^{-{1\over a}}) {X_1}(1-{\lambda}^{-{1\over a}})^{{X_1}-1}\big).$$
When $\lambda<2^a$, $p=1/2$ and $f_1,f_2$ are defined by
$$f_1(x,y,z)=-{{xz}\over{(1-x)y\log y}},\qquad f_2(x,y,z)=-x^{-z}\log y.$$
In this case, $W_1$ and $W_2$ are given by
$$W_1=\lambda^{-1}2^a
e^{{\lambda}\over{2^a}}(X_12^{-X_1}+a(1-\lambda 2^{-a})2^{-X_1})$$
and $$W_2=2^a
e^{{\lambda}\over{2^a}}(X_12^{-X_1}\log 2+\big(a(1-\lambda 2^{-a})\log 2-1\big)2^{-X_1}).$$
Corollary \ref{cor} is so proven.
\qed

\end{document}